\documentclass[12pt,reqno]{amsart}

\usepackage{latexsym}

\newcommand{\T}{{\mathbf T}^m}

\newcommand{\sm}{\setminus}
\newcommand{\szego}{Szeg\"o }

\newcommand{\Si}{\Sigma}
\newcommand{\inv}{^{-1}}
\newcommand{\kahler}{K\"ahler }

\newcommand{\wt}{\widetilde}
\newcommand{\wh}{\widehat}
\newcommand{\PP}{{\mathbb P}}
\newcommand{\N}{{\mathbb N}}
\newcommand{\R}{{\mathbb R}}
\newcommand{\C}{{\mathbb C}}

\newcommand{\Z}{{\mathbb Z}}

\newcommand{\CP}{\C\PP}
\renewcommand{\d}{\partial}
\newcommand{\dbar}{\bar\partial}
\newcommand{\ddbar}{\partial\dbar}

\newcommand{\half}{{\frac{1}{2}}}
\newcommand{\vol}{{\operatorname{Vol}}}

\newcommand{\FS}{{{\operatorname{FS}}}}

\renewcommand{\phi}{\varphi}

\newcommand{\ccal}{\mathcal{C}}

\newcommand{\hcal}{\mathcal{H}}

\newcommand{\lcal}{\mathcal{L}}

\newcommand{\ocal}{\mathcal{O}}
\newcommand{\pcal}{\mathcal{P}}
\newcommand{\qcal}{\mathcal{Q}}

\newcommand{\al}{\alpha}
\newcommand{\be}{\beta}
\newcommand{\ga}{\gamma}

\newcommand{\La}{\Lambda}
\newcommand{\la}{\lambda}

\newcommand{\om}{\omega}

\newtheorem{theo}{{\sc Theorem}}[section]

\newtheorem{lem}[theo]{{\sc Lemma}}
\newtheorem{prop}[theo]{{\sc Proposition}}

\newenvironment{rem}{\medskip\noindent{\it Remark:\/} }{\medskip}
\newenvironment{defin}{\medskip\noindent{\it Definition:\/} }{\medskip}

\title[Harmonic Analysis on Toric Varieties] {Harmonic Analysis on Toric
Varieties}

\author{Bernard Shiffman}

\address{Department of Mathematics, Johns Hopkins University, Baltimore, MD
21218, USA} \email{shiffman@math.jhu.edu}

\author{Tatsuya Tate}
\address{Department of Mathematics, Keio University, Keio University
3-14-1 Hiyoshi Kohoku-ku, Yokohama, 223--8522 Japan}
\email{tate@math.keio.ac.jp}

\author{Steve Zelditch}
\address{Department of Mathematics, Johns Hopkins University, Baltimore, MD
21218, USA} \email{zelditch@math.jhu.edu}

\thanks{Research partially supported by NSF grants DMS-0100474 (first
author) and  DMS-0071358 (third author) and by JSPS (second
author).}

\date{February 28, 2003}

\begin{document}

\begin{abstract} Harmonic analysis on a toric \kahler variety $(M, \omega)$
refers to the orthonormal basis  of  eigenfunctions of the
$(\C^*)^m$  action on
 the spaces $H^0(M, L^N)$ of holomorphic sections of powers of the positive
line bundle
  $L \to M$ with $c_1(L) = [\omega]$,  and the Fourier multipliers that
act on them.
 Using this harmonic analysis, we give an   exact formula
  for the \szego kernel  as a
Fourier multiplier ${\mathcal M}$ applied to the pull back  of the
\szego kernel of projective space under a
 monomial embedding. The Fourier multiplier ${\mathcal M}$ involves a
partition function of the convex lattice polytope $P$ associated
to $M$. We further prove that ${\mathcal M}$ is  a Toeplitz
operator, and as a corollary we obtain an oscillatory integral
formula for the characters $\chi_{NP}$ of the torus action on
$H^0(M, L^N)$.

\end{abstract}

\maketitle

\section{Introduction}
Toric \kahler varieties $(M, \omega)$ are often used in geometry
as simple test cases  for  difficult geometric problems, e.g.
mirror symmetry or existence of K\"ahler-Einstein metrics.  Their
simplifying features are the existence of a holomorphic action of
$(\C^*)^m$ ($m  = \dim M)$ with an open dense orbit,  and of a
Hamiltonian action of the underlying  real torus $\T$  with
respect to the \kahler form. The associated  moment map $\mu: M
\to P \subset \R^m$ expresses
 $M$ as a stratified torus fibration over the
convex lattice polytope $P$. In short, toric varieties are the
completely integrable systems of complex analysis and many
interesting quantities are explicitly solvable on them, often in
terms of the combinatorics of $P$.

 The purpose of this article is to
consider  toric varieties as a model setting for harmonic analysis
on \kahler manifolds.  Motivated by complex analysis, we consider
a Hermitian holomorphic line bundle $L \to M$ with $c_1(L) =
\omega$. Harmonic analysis in this paper refers to the Hilbert
space completion of the coordinate ring,
\begin{equation} {\mathcal H}:= \bigoplus_{N = 0}^\infty\ H^0(M,
L^N)\;, \end{equation} where $H^0(M, L^N)$ denotes the space of
holomorphic sections of the $N$-th  tensor power of $L$. The torus
action gives rise to a natural Fourier analysis on ${\mathcal H}$.
In fact, it extends to all $\lcal^2$ sections.

As in the classical settings of Fourier analysis,
square-integrable functions on $\R^n$ and on the real $n$-torus
${\bf T}^n = \R^n/\Z^n$, the space ${\mathcal H}$  is spanned by
exponentials, or more precisely, eigenfunctions of the linearized
$(\C^*)^m$ action. In other coordinates, they appear as
`monomials.' In the fundamental case of $M = \CP^m$, the joint
eigenfunctions are  the monomials given in an affine chart by
\begin{equation} \chi_{\alpha} : \C^m \to \C,\;\;\;
\chi_{\alpha}(z) = z^{\alpha}. \end{equation} The monomials lift
(by homogenization) to homogeneous monomials on $\C^{m + 1}$. In
general, the linearized ${\bf T}^m$ action is generated by  $m$
commuting operators $\hat{I}_j$, $j = 1, \dots, m$ on $M$ which
preserve $H^0(M, L^N)$, and the joint spectrum of the eigenvalue
problem
\begin{equation} \label{QCI} \hat{I}_j \phi_{\alpha} = \alpha_j
\phi_{\alpha},\;\; \alpha \in \R^m,\;\;\; \bar{\partial}
\phi_{\alpha} = 0
\end{equation}  consists of lattice points $\alpha \in N P \cap \Z^m$.
All of $\lcal^2(M,L^N)$ is spanned by (not necessarily
holomorphic) monomials, but we will only be studying the
holomorphic ones in detail.

Toric varieties are thus models of quantum completely integrable
systems. They are of a very special type because integrable
systems usually generate an $\R^m$ action rather than a ${\bf
T}^m$ action.  On such varieties, one can obtain analytical
results with more precision than is possible in almost any other
case. In this paper, we will illustrate this theme with regard to
the \szego kernels $\Pi_N=\Pi_N^M: \lcal^2(M, L^N) \to H^0(M, L^N)$.

Our first result is an exact formula for the  \szego kernel of a
smooth, projective,  toric \kahler variety.  We shall define our toric varieties
through  monomial embeddings as follows (see \cite[Chapter~5]{GKZ}): Let
$P$ be a convex integral polytope in $\R^m$ and denote the lattice points in $P$  by
$P\cap\Z^m=\{\al(1),\al(2),\dots\al(\#P)\}$. For simplicity, we shall assume
throughout this paper that $P$ is contained in the positive quadrant
$[0,+\infty)^m$. For any vector
$c = (c_1, \dots, c_{\# P})\in(\C^*)^{\#P}$, we define the  map
\begin{equation} \label{PHIPC} \Phi_P^c = \big[c_{\al(1)}\chi_{\al(1)},\dots
,c_{\al(\#P)}\chi_{\al(\#P)}\big]:(\C^*)^m\to\CP^{\#P-1} \;;
\end{equation}  i.e.,
$$ \Phi_P^c(z)= \big[c_{\al(1)}z^{\al(1)},\dots
,c_{\al(\#P)}z^{\al(\#P)}\big]\;.$$ The closure of the image is the
toric variety $M_P^c \subset \CP^{\# P - 1}$.  If $P$ is Delzant,
then $\Phi_P^c$ is an embedding,  and we identify $\C^{*m}$ with
its image (the `open orbit') in $M_P^c$. To simplify the notation
when $c$ and $P$ are fixed,  we  often write the  resulting
embedding as
\begin{equation}  \Phi: M \hookrightarrow \CP^{d},\;\;
d=\#P -1.
\end{equation} and refer to it as  a monomial embedding. We further define
the line bundle
\begin{equation} L =  L_P^c := \Phi^* {\mathcal O}(1).
\label{LP}\end{equation}

Our analysis is based on identifying $\hcal$ with the Hardy
space of CR functions,
$$\hcal^2(X):=\{F\in\lcal^2(X): \dbar_bF=0\}\;,$$ where $X$ is the unit circle
bundle in the dual line bundle $L\inv$. Its $N$-th Fourier component
$$\hcal^2_N(X)= \{F\in
\hcal^2(X): F(e^{i\theta}x)=e^{iN\theta}F(x)\}\;,$$ can be identified with
$H^0(M,L^N)$ via the equivariant lifting of holomorphic sections to $X$ (see \S
\ref{FAnalysis}). The monomial embedding $\Phi$ likewise lifts to an embedding of the
$S^{1}$ bundle $X$:
\begin{equation}\label{TAM}
\iota:X \hookrightarrow S^{2d+1}.
\end{equation}
(See \S \ref{s-lifting}.)

Our  explicit formula for the \szego kernels $\Pi_N$ of the line
bundles $L^N \to M$, i.e.\ the orthogonal projections onto the
spaces $H^0(M, L^N)$, involves two ingredients. The first is the
pullback of the projective space \szego kernels of powers ${\mathcal O}(N) \to
\CP^{\#P - 1}$ under the monomial embedding:
\begin{equation} \iota^* \Pi_N^{\CP^{d}} (x, y)
= \frac{(N+m)!}{N!}\, \langle \iota(x), \overline{\iota(y) }\rangle^N,\qquad
x,y\in X\; .
\label{iota}\end{equation} As we will see, these simple kernels  are the key
objects in all of our results.

Our first result on the relation of $\Pi=\sum_{N=0}^\infty\Pi_N^M$ to $\Pi_{\iota} =
\iota^* \Pi^{\CP^d}$ is valid on all integral \kahler manifolds.
In Lemma \ref{PIOTA}, we prove that there always exists a
pseudo-differential operator $A$ of order $m$ such that $\Pi \sim
\Pi  A  \Pi \Pi_{\iota} = \Pi  A \Pi_{\iota}  $ modulo smoothing
operators $\Pi R \Pi.$ Note that this  relation is equivalent to
$\Pi \sim A \Pi_{\iota}$ since we can left multiply the latter by
$\Pi$. Thus, only the Toeplitz part of $A$, i.e. its restriction
to the Hardy space, is important in this relation. This relation
might be contrasted with the comparison inequality  of Li-Tian  \cite{LT},
which says that the heat kernel of a \kahler manifold is less than
the pull back of the projective space heat kernel under an
isometric holomorphic embedding. When $t \to \infty$
this implies an inequality between $\Pi$ and $\Pi_{\iota}$.

 The two kernels are  more  closely related on a smooth
toric variety because the sections of $L_P^N$ are spanned by
products of the monomials $\chi_\al$ for $\al\in P\cap \Z^m$. However, $\Pi_N$
is not equal to $\iota^* \Pi_N^{\CP^d}$ and they even have
different orders as complex Fourier integral operators. But our
next result shows that one may  adjust one into the other using
{\it Fourier multipliers} (or convolution operators)  on a toric
variety. The definition is analogous to the Euclidean one, namely
an operator which commutes with the group action,  or equivalently
has the exponential functions as its eigenfunctions.

Our main results are summarized in the following theorem:

\begin{theo}\label{M} Let $(M_{P}, \omega_{P})$ be a \kahler toric variety.
Then there   exists a Fourier multiplier ${\mathcal M}$ such that
$$\Pi_N^{M_P} (x, y) = {\mathcal M} \; \cdot\;  \langle \iota(x), \overline{\iota(y)}
\rangle^N, $$ where $\iota: X_P \to S^{2d+1}$ is the lift of the
monomial embedding $\Phi$ for which $\Phi^* \omega_{FS} =
\omega_{P}.$ Moreover, $\Pi^{M_P} {\mathcal M}\, \Pi^{M_P}$ is a Toeplitz
operator.

\end{theo}

 The Fourier multiplier is
defined by the condition that  its eigenvalue on a joint
eigenfunction $\phi_{\alpha} \in H^0(M_P, L_P^N)$ with $\alpha \in
N P$ is given by \begin{equation} \label{FM} {\mathcal
P}_N^{-1}(\alpha) {\mathcal Q}_N^{-1} (\alpha) \end{equation}
where:

\begin{enumerate}

\item [(i)]  ${\mathcal P}_N$ is the lattice path  `partition function'
$${\mathcal P}_N(\alpha)  = \# \{(\beta_1,
\dots, \beta_N): \beta_j \in P, \beta_1 + \cdots + \beta_N =
\alpha \}. $$

\item [(ii)] ${\mathcal
Q}_N(\alpha)$  is the monomial norming function: $$ {\mathcal
Q}_N(\alpha):= \int_{M_P} \|\chi^P_\al (z)\|_{h_P^N}^2 \,d\vol_{M_P}(z).$$
(Here, $\chi^P_\al$ is the section of $L_P^N$ corresponding to $\chi_\al$; see \S
\ref{s-bundles}.)

\end{enumerate}

The statement that $\Pi {\mathcal M} \Pi$ is a Toeplitz operator
means that there is a smooth symbol $\sigma $ in the usual
semiclassical sense so that $\Pi {\mathcal M} \Pi = \Pi \sigma
\Pi$, modulo smoothing operators. Thus, the multiplier has simple
asymptotic properties, e.g. it is polyhomogeneous  along rays of
lattice points. The lattice path partition function seems to be
of some interest in its own right. In a subsequent paper, we give
asymptotic formulae for lattice path counting functions and
applications to multiplicities of group representations
\cite{STZ}.

As a consequence of our proof of Theorem \ref{M}, we  obtain  the following integral
formula for the equivariant characters $\chi_{NP}$ on $(\C^*)^m$ given by
\begin{equation} \chi_{N P} (e^{i\phi}) = \sum_{\alpha \in N P}
e^{i\langle \phi,\al\rangle}\ ,\qquad \phi=(\phi_1,\dots,\phi_m).
\label{char}\end{equation}

\begin{theo} The  characters $\{\chi_{NP}\}$ are  given by  a  complex
oscillatory integral of the form $$\chi_{NP} (e^{i\phi}) \sim \int_{M_P} e^{N \Psi(x,
\phi)} A_N(x) dV(x), $$ where $\Psi$ is a non-degenerate
complex phase function of positive type and  $A_N$ is a
polyhomogeneous function of $N$.
\end{theo}

Here,  `$\sim$' means modulo a smoothing operator. A more precise formula is
given as Proposition \ref{aha}. We note that another oscillatory integral formula
for the polytope character was given in
\cite[Prop.~2.1]{SZ2} using completely different methods.  In the formula in
\cite{SZ2}, the amplitude is a polynomial in $N$, but depends on $\theta$ as well as
$x$.  Proposition~2.1 in \cite{SZ2} was used  to obtain asymptotic formulas for the
distribution of zeros of random polynomials with expanding Newton polytope $NP$. In
fact, Proposition
\ref{aha} provides an alternate approach to the results on zero distributions in
\cite{SZ2} (for the case of Delzant polytopes) as well as to the results in our
forthcoming paper
\cite{STZ1} on the distribution of values of eigenfunctions on toric varieties.

\medskip
\noindent{{\it Acknowledgments:\/}} We would like to thank Amit Khetan for helpful
comments concerning projective normality for toric varieties. This paper was written
during a stay of the second author at Johns Hopkins University. He would like to
express his special thanks to the faculty in the Department of Mathematics of Johns
Hopkins University.

\medskip

\section{Background on toric varieties and moment polytopes}

 Recall that a {\it toric variety\/} is a complex
algebraic variety $M$ containing the {\it complex torus\/} $$\C^{*m}:=
(\C\sm\{0\})\times\cdots\times (\C\sm\{0\})$$ as a Zariski-dense
open set such that the group action of $\C^{*m}$
 on itself extends to a $\C^{*m}$ action on $M$. In the smooth case, $M$ can be
given the structure of a symplectic manifold such that the
restriction of the action to the underlying real torus
$$\T=\{(\zeta_1,\dots\,\zeta_m)\in \C^{*m}: |\zeta_j|=1, 1\le j\le m\}$$ is a
Hamiltonian action (see \S \ref{s-torus}). As mentioned in the
introduction, we  define projective toric varieties  by the monomial
embeddings (\ref{PHIPC}), which depend on a convex lattice
polytope $P$ and a choice of weights $c_{\alpha}$. In the smooth
case, $P$ satisfies Delzant's condition (see \cite{Gu} or \cite{SZ2} for the
definition).

\subsection{The bundles $L_P^c$ and their holomorphic sections}
\label{s-bundles}

Recall that $H^0(\CP^{\#P-1},\ocal(1))$ has as a basis the linear
coordinate functions $\la_j:\C^{\#P}\to\C$, $1\le j\le {\#P}$. The
{\it Fubini-study metric\/} $h_\FS$ on $\ocal(1)$ is given by
$$|\la_j|_\FS ([\zeta])=\frac{|\zeta_j| }{\|\zeta\|}\qquad
(\zeta\in\C^{\#P})\;,$$ which has curvature form
$\om_\FS=\frac{i}{2\pi}\ddbar \log \|\zeta\|^2$. We  endow $L_P^c$
with the Hermitian metric
 $h_P^c:= \Phi_P^{c*} h_{\FS}$ of curvature $\omega_P^c$ given on
$\C^{*m}$ by
\begin{equation} \om_P^c=\Phi_P^{c*}\om_\FS=\frac{\sqrt{-1}}{2\pi}\partial
\bar{\partial} \log \sum_{\alpha \in P} |c_{\alpha}|^2
|z^{\alpha}|^2. \label{pullfubini}
\end{equation}

Each monomial $\chi_{\alpha}$  with $\alpha \in P$ corresponds to
a section of $H^0(M_P^c, L_P^c)$ and vice versa. To explicitly
define this correspondence, we make the identifications (recalling
(\ref{LP})):
\begin{equation}\label{identify}\chi^{P}_{\al(j)}\equiv
c_{\al(j)}\inv\Phi_P^c{}^*\zeta_j\in H^0(M_P^c, L_P^c)\;, \quad
1\le j \le \#P\;.\end{equation}
So far, we have not specified the constants $c_\al$.  For studying
our phenomena, the choice of constants defining the toric variety
$M_P$ is not important.   (We use $c_\al= 1$ in \S
\ref{SKTV}, but except for Lemma \ref{PQ},  any choice of the $c_\al$ will work.)

More generally, a basis for the space $H^0(M_P,L_P^N)$ of global sections
of
$L_P^N$ is given by the monomials  $\chi^{P}_{\gamma}$, where $\ga$ runs over the
lattice points of
$NP$. See, for example, \cite[\S
3.4]{Fu}.  (In \cite{Fu}, the toric variety given by $P$ is defined using the
normal fan of $P$.  However  by a theorem of Demazure \cite[p.~71]{Fu}, the
associated line bundle in this construction is very ample since $M_P$ is smooth, and
hence the fan and monomial embedding constructions are algebraically equivalent.)

We note that since $P$ is Delzant, each $\ga\in NP\cap\Z^m$ can be decomposed as
\begin{equation}\label{solid}
\chi^{P}_{\gamma}=\chi^{P}_{\beta_{1}} \otimes \cdots \otimes \chi^{P}_{\beta_{N}},
\end{equation}
where $\beta_{1},\ldots ,\beta_{N} \in P\cap\Z^m$ such that $\gamma
=\beta_{1}+ \cdots +\beta_{N}$. Such a partition of $\gamma$ exists since $M_P$ is
smooth and hence projectively normal \cite[pp.~72--73]{Fu}; i.e., the cone
$C_P\subset
\C^{\#P}$ over
$M_P$ is normal. Indeed, by normality of $C_P$, each section of $H^0(M_P,L_P^N)$
corresponds to a homogeneous polynomial of order $N$ on   $\C^{\#P}$; i.e.,
$H^0(M_P,L_P)$ generates the coordinate ring $\bigoplus_{N=1}^\infty
H^0(M_P,L_P^N)$, and (\ref{solid}) follows.

\begin{rem} For a general integral polytope $P$, (\ref{solid}) may not hold.
A well-known example is where $P$ is the simplex in $\R^3$ with vertices $(0,0,0),\
(0,1,1),\ (1,0,1),\ (1,1,0)$.   Then $P\cap\Z^3$ consists only of the vertices, and
for all $N\ge 2$, the lattice point $\ga=(1,1,1)$ lies in $NP$ but cannot be
decomposed as in (\ref{solid}).\end{rem}

\subsection{Lifting to the associated $S^1$ bundle}\label{s-lifting}

The geometry of line bundles can be rephrased in terms of the
 associated principal $S^1$ bundle
$$X_P  = \{v \in L_P\inv: \|v\|_{h_P\inv} =1\},$$
where $L_P\inv \to M_P$ denote the dual line bundle to $L_P$ with dual
metric $h_P\inv$.

\medskip The action of the real torus ${\bf T}^m$ lifts from
$M_P^c$ to $X_P^c$ and combines with the $S^1$ action to define a
${\mathbf T}^{m+1}$ action on $X_P^c$.  Recall that under the
monomial embedding
$$\Phi_P^c:\C^{*m}\hookrightarrow M_P^c\hookrightarrow \CP^{\#P-1},\qquad
z\mapsto \big[c_{\al(1)}z^{\al(1)},\dots
,c_{\al(\#P)}z^{\al(\#P)}\big]\;,$$ the $\T$ action on $M_P^c
\subset \CP^{\#P-1}$ is given by
\begin{equation}\label{actionM}e^{i\phi} \cdot
[\zeta_1,\dots,\zeta_{\#P}]=
\big[e^{i\langle\al(1),\phi\rangle}\zeta_1,\dots,
e^{i\langle\al(\#P),\phi\rangle}\zeta_{\#P}\big]\;, \quad
e^{i\phi}=(e^{i\phi_1}, \dots, e^{i\phi_m}).\end{equation} The action (\ref{actionM})
lifts to an action  on
$L_P\inv$:
\begin{equation}\label{actionX}e^{i\phi} \cdot \zeta = \big( e^{i\langle\al(1),\phi\rangle}\zeta_1,\dots,
e^{i\langle\al(\#P),\phi\rangle}\zeta_{\#P}\big) \;.\end{equation}
Since the circle bundle $X_P^c\subset S^{2\#P-1}$ is invariant
under this action, (\ref{actionX}) also gives a lift of the action
(\ref{actionM}) to $X_P^c$.

We also have the standard circle action on $X_P^c$:
\begin{equation}\label{circle} e^{i\theta}\cdot \zeta =  e^{i\theta}
\zeta\;,\end{equation} which commutes with the $\T$-action
(\ref{actionX}). Combining (\ref{actionX}) and (\ref{circle}), we
then obtain a ${\mathbf T}^{m+1}$-action on $X_P^c$:
\begin{equation}\label{actionX+}(e^{i\theta}, e^{i\phi_1},\dots, e^{i\phi_m})\bullet
\zeta =  e^{i\theta}(e^{i\phi}\cdot \zeta)\;.\end{equation}

We
identify sections $s_N$ of $L^N$ with equivariant functions
$\hat{s}_N$ on $X$ by the rule
\begin{equation} \label{sNhat}\hat{s}_N(\lambda) = \left( \lambda^{\otimes
N}, s_N(z) \right)\,,\quad \la\in X_z\,,\end{equation} where
$\lambda^{\otimes N} = \lambda \otimes \cdots\otimes \lambda$.
Clearly,   $\hat s_N(e^{i \theta} \cdot  x) = e^{iN\theta} \hat
s_N(x)$ if  $s \in H^0(M_P^c, L_P^N)$.

It should be noted that, for each $s \in H^{0}(M_{P},L_{P}^{N})$, we have
\begin{equation}
|\wh{s}_{N}(x)|=\|s(z)\|_{h_P^N}
\;.\end{equation}

We now introduce notation for the lifts of monomials: for $\al\in P$, we lift
$\chi_{\al}^{P}\in H^0(M_P^c, L_P^c)$ to an equivariant function
$\wh\chi_\al^P$ on the circle bundle $X_P^c\to M_P^c$, and we
write
\begin{equation}\label{mhat}\wh m_{\al(j)}^P:= c_{\al(j)} \wh
\chi_{\al(j)}^P= \zeta_j\circ\iota_P\end{equation} where $\iota_P
:X_P^c \to S^{2d +1}$ is the lift of the  embedding  $M_P^c
\hookrightarrow \CP^{d}$ ($d=\#P-1$). (Of course, $\wh m_\al^P$ depends on
$c$, which we omit to simplify notation.) We also consider the
monomials
$$m_{\al}^{P}:=c_\al \chi_{\al}^{P}$$
so that $\wh m_\al^P$ is the equivariant lift of $m_\al$ to
$X_P^c$.  In terms of local coordinates $(z,\theta)$ on
$\pi\inv(\C^{*m})\subset X_P^c$, we have
\begin{equation}\label{equivariantm}
\wh m_\al^P(z,\theta)= \frac{e^{i\theta}c_\al
z^\al}{\left(\sum_{\be\in P}|c_\be|^2|z^\be|^2\right)^{1/2}}\;.
\end{equation}

\subsection{Moment maps and torus actions}\label{s-torus}

The group $\C^{*m}$ acts on $M_P^c$ and the subgroup $\T $ acts in
a Hamiltonian fashion. Let us recall the formula for its moment
map $\mu_{P}^c: M_P^c \to \R^m$,  restricted to the open orbit
$\C^{*m}$. This  moment map is the composition
$$\mu_P^c: M_P^c\subset \CP^{\#P -1} \buildrel{\mu_0}\over \to
\R^{\#P } \buildrel{A}\over\to \R^m\;,$$ where
$$\mu_0([z_1, \dots, z_{\#P}]) = \frac{1}{\|z\|^2} (|z_1|^2, \dots,
|z_{\#P}|^2),$$ and $A$ is the linear projection is given by the
column vectors $(\alpha^1, \dots, \alpha^{\#P})$.  Hence we have
\begin{equation}\label{muP} \mu_P^c = \frac{1}{\sum_{\alpha \in P}
|c_{\alpha}|^2 |\chi_{\alpha}|^{2 }} \sum_{\alpha \in P}
|c_{\alpha}|^2 |\chi_{\alpha}|^{2 }\alpha \;. \end{equation} For
any $c$, the image of $M^c_P$ under $\mu_P^c$ equals $P$.

Noting that
\begin{equation}\label{sum=1}\sum_{\al\in P}|\wh m_\al^P|^2 =\sum_{j=1}^{\#P}|\zeta_j\circ\iota_P|^2 \equiv 1\;,\end{equation}
we obtain the formula:
\begin{equation}\label{muP2} \mu_P^c(z) = \sum_{\alpha \in P}
|\wh m_\al^P(z)|^2\alpha\;.\end{equation} (We write $|\wh
m_\al^P(z)|=|\wh m_\al^P(z,\theta)|$, since the absolute value is
independent of $\theta$.)

\subsubsection{Projective space} We illustrate with the case of
$\CP^m$ and introduce
notation that will be used throughout the paper. As a toric variety, $\CP^m=M_\Si$,
where
$\Si$ is the standard simplex in $\R^m$ with vertices at the points
$$(0,\dots,0),\ (1,0,\dots, 0),\ (0,1,\dots,0),\ \dots,(0,\dots,0,1)\ ,$$
and $L_\Si$ is
the hyperplane section bundle $\ocal(1)\to\CP^m$. In this case, $X
= S^{2m + 1}$ and the lifts of sections in $H^0(\CP^m,\ocal(p))$
to $X$  consist of homogeneous polynomials
$$F(\zeta_0, \dots,
\zeta_m) = \sum_{|\la| = p} C_\la \zeta^\la \qquad (\zeta^\la =
\zeta_0^{\la_0} \cdots \zeta_m^{\la_m})$$ in $m+1$ variables. We
give $\CP^m$ the Fubini-Study \kahler form given in homogeneous
coordinates $(\zeta_0,\dots,\zeta_m)$ by $\om_\FS=
\frac{i}{2\pi}\ddbar \log\|\zeta\|^2$, and we   give $\ocal(p)$
the {\it Fubini-Study metric\/}:
$$|F(\zeta)|_\FS= |F(\zeta)|/\|\zeta\|^p, \quad\mbox{for }\ F\in
H^0(\CP^m,\ocal(p))\;.$$ Identifying  $F$ with the polynomial
$f(z)=F(1,z_1,\dots,z_m)$, the Fubini-Study norm can be written
$$|f(z)|_\FS=|f(z)|/(1+\|z\|^2)^{p/2} \qquad (z\in \C^m)\;.$$

We equip the space $H^0(\CP^m, \ocal(p))$ of all
homogeneous polynomials of degree $p$ with the Hermitian inner product:
\begin{equation}\label{IP}\langle f, \bar g \rangle = \int _{\CP^m}\left\langle
F,\overline{G}\right\rangle_\FS \,d\vol_{\CP^m}=  \frac{1}{m!}
\int _{\C^m}\frac{\langle
f(z),\overline{g(z)}\rangle}{(1+\|z\|^2)^p} \,\om_\FS^m(z),\quad
f,g\in  H^0(\CP^m, \ocal(p)).\end{equation} (We use here the
Riemannian volume $d\vol_{\CP^m}=\frac{1}{m!}\om_\FS^m$; note that
the total volume of $\CP^m$ is $\frac{1}{m!}$, using our
conventions.)

Under the $\T $ action, we have the weight space decomposition
 $$H^0(\CP^m, {\mathcal O}(p)) = \bigoplus_{|\alpha|\le p}  \C
\chi_{\alpha}\;,$$ where we recall that $\chi_\al(z) =
z_1^{\alpha_1} \cdots z_m^{\alpha_m}$. The monomials  $\{
\chi_{\alpha}\}$ are orthogonal but not normalized. Any choice of
norming constants $\{r_{\alpha} \in \C^{*}\}$ will give a monomial
basis $\{  r_{\alpha} \chi_{\alpha}\}$ for $H^0(\CP^m, {\mathcal
O}(p))$. We shall choose $r_\al=\|\chi_\al\|_{\CP^m}\inv$,  where
$$\|\chi_\al\|_{\CP^m} =\sqrt{ \langle \chi_\al, \chi_\al\rangle} =
\left[\frac{p!}{(p+m)!{p\choose\al}}\right]^\half\;,\quad\quad
{p\choose\al}:= \frac{p!}{(p-|\al|)!\alpha_1!\cdots \alpha_m!}\;,
$$ is the Fubini-Study $\lcal^2$ norm of $\chi_{\alpha}$
given by (\ref{IP}). (See \cite[\S 4.2]{SZ1}; the extra factor $m!$
in \cite{SZ1} is due to the use of $\om^m$ instead of
$\frac{\om^m}{m!}$ for the volume form.) This choice provides an
orthonormal basis for $H^0(\CP^m, {\mathcal O}(p))$ given by the
monomials
$$\frac{1}{\|\chi_\al\|}_{\CP^m}\,\chi_\al= \sqrt{\frac{(p+m)!}{p!} {p\choose
\al}}\ \chi_\al\ ,\qquad |\al|\le p\;.$$

In this case, we shall use the special choice of the coefficients
of the monomial embedding
$$c^*_\al:= {p\choose
\al}^\half=\left(\frac{p!}{(p+m)!}\right)^\half\;\|\chi_\al\|_{\CP^m}\inv\;,$$
so that
\begin{equation}\label{special} \mu_{p\Si}(z):=\mu_{p\Si}^{c^*}
(z)=  \frac{1}{\sum_{|\alpha|\le p} {p\choose\al}|z^{\alpha}|^{2
}} \sum_{|\alpha|\le p} {p\choose\al}|z^{\alpha}|^{2 }\alpha
=\frac{p}{1+\sum|z_j|^2}(|z_1|^2,\dots ,|z_m|^2)\,,\end{equation}
where the last equality follows by differentiating the identity
$(1+\sum x_j)^p=\sum_{|\al|\le p} {p\choose\al}x^\al$.  Note that
this choice gives us the scaling formula
$$\mu_{p\Si}=p\mu_\Si.$$

We can identify $L\inv=\ocal_{\CP^{m}}(-1)$ with $\mathbb{C}^{m+1}$ with the
origin blown up, and the circle bundle $X \subset L\inv$ is identified with
the unit sphere $S^{2m +1} \subset \mathbb{C}^{m+1}$. Then the equivariant lift
$\wh{\chi}_{\alpha}:S^{2m+1} \to \mathbb{C}$ of $\chi_{\alpha} \in H^{0}(\CP^{m},\ocal(p))$
is given by the homogenization:
\[
\wh{\chi}_{\alpha}(x)=x^{\wh{\alpha}},\quad
\wh{\alpha}=(\alpha_{1},\ldots,\alpha_{m},p-|\alpha|).
\]

Furthermore, \begin{equation}\label{mchi}|\wh
m_\al^{p\Si}(z)|=\left[ \frac{p!}{(m+p)!}\right]^\half \frac{|\wh
\chi_\al(z)|}{\ \|\chi_\al\|_{\CP^m}}\;,\end{equation} where, by abuse of
notation, we regard $|\wh\chi_\al|$ and $|\wh m_\al^{p\Si}|$ as
functions on $\C^{*m}$, since they are invariant under the circle
action.
Strictly speaking, the $S^{1}$-bundle $X_{p\Si}$ in this case is given by
$X_{p\Si}\equiv S^{2m+1}/\Z_p$. Then we have the
more precise formula
$$\wh m_\al^{p\Si}(x')=\left[
\frac{p!}{(m+p)!}\right]^\half\frac{\wh
\chi_\al(x)}{\ \|\chi_\al\|_{\CP^m}}\;,$$ where $x'\in M_{p\Si}$ is the
equivalence class of $x\in S^{2m+1}$,  and therefore we have
\begin{equation}\label{msum}\sum_{|\al|\le p}\wh m_\al^{p\Si}(x') \overline {\wh m_\al^{p\Si}(y')} =\langle x,\bar
y\rangle^p \;.\end{equation}

\section{Fourier analysis}
\label{FAnalysis}

As mentioned in the introduction, we
identify ${\mathcal H}$ with the Hardy space $\hcal^2(X_P^c)$. It
is a Hilbert space with the inner product $$\langle f, \bar g \rangle =
\int_X f\bar g \,d\vol_X,\;\;\; d\vol_X = \frac 1{m!}\alpha \wedge (d\alpha)^m, $$
where $\alpha$ is a contact 1-form defined by the Hermitian connection on $L_{P}\inv$
such that $d\alpha =\pi^{*}\omega_{P}$.
Under the identification of sections with equivariant functions,
the inner product is the same as
\begin{equation}\langle s_1,\bar s_2 \rangle_{N} = \int_M \left\langle
s_1(z),\overline{s_2
(z)}\right\rangle_{h_P^N}d\vol_M(z)\;,\quad\quad s_1,s_2\in
H^{0}(M_{P},L_{P}^{N})\;.\label{inner}\end{equation}

 Under the $S^1$ action, the Hardy space has the orthogonal
decomposition
\begin{equation} \hcal^2(X_P^c)=
\bigoplus_{N=0}^{\infty}\hcal^2_N(X_P^c). \end{equation} We recall
that the (equivariant) `\szego projectors' $\Pi_{N}$ are the
orthogonal projection onto $H^0(M_P^c, L_P^{c N})$. If $\{S^N_j\}$
denotes an orthonormal basis of this space, and $\hat{S_j}^N$
denote their lifts to $X$,  then the projector $\Pi_{N}$ is given
by the kernel
\begin{equation}\label{szego}\Pi_{N}(x,y) = \sum_{j = 1}^{k_N} \hat
S_j^N(z) \overline{\hat S_j^N(y)}\;: \lcal^2(X_P^c) \to
\hcal^2_N(X_P^c).\end{equation} The full \szego kernel is the
equivariant direct sum \begin{equation} \label{Xszego} \Pi =
\sum_{N = 1}^{\infty} \Pi_N:  \lcal^2(X_P^c) \to
\hcal^2(X_P^c).\end{equation}

For each $N \in \mathbb{N}$, we define the `homogenization'
$\wh{NP} \subset \mathbb{Z}^{m+1}$
of the lattice point in the polytope $NP$ to be the set of all lattice point
$\wh{\alpha}^{N}$ of the form
\[
\wh{\alpha}^{N}=\wh{\alpha}:=(\alpha_{1},\ldots,\alpha_{m},Np -|\alpha|),\quad
\alpha=(\alpha_{1},\ldots,\alpha_{m}) \in NP \cap \mathbb{Z}^{m},
\]
where, as before, we set $p=\max_{\beta \in P \cap \mathbb{Z}^{m}}|\beta|$.
We also define the cone $\Lambda_{P}=\bigcup_{N=1}^{\infty}\wh{NP}$.

In this section, we consider the lifted monomials by $\wh{\chi}_{\wh{\alpha}}
(x)=\wh{\chi}_{\alpha}
(x)$, for $\wh{\alpha} \in \Lambda_{P}$.
We may combine the problems as $N$ varies into a homogeneous
eigenvalue problem on $X$:  After homogenization, i.e. lifting to
$X$, we obtain the joint scalar eigenvalue problem
\begin{equation} \label{jointeigen}
\hat{I}_j \wh{\chi}_{\wh{\alpha}} = \wh{\alpha}_j \wh{\chi}_{\wh{\alpha}},\;\;
\wh{\alpha} \in \R^{m + 1} ,\;\;\; \bar{\partial}_b
\wh{\chi}_{\wh{\alpha}} = 0, \;\;
j=1,\ldots,m+1.
\end{equation}

  The lattice points $\wh{\alpha}$ lie in the cone $\Lambda_P \subset \R^{m
+ 1}$. It is well known that rays $\N \wh{\alpha}$ in this cone define
a semiclassical limit.

The torus action on $X_P^c$ can be quantized to define an action
of the torus as unitary operators on $\hcal^2(X_P^c).$
Specifically, we let $\Xi_1,\dots\Xi_m$ denote the differential
operators on $X_P^c$ generated by the $\T$ action:
\begin{equation}\label{Xi}(\Xi_j \hat S)
(\zeta)=\frac{1}{i}\frac{\d}{\d\phi_j}  \hat
S(e^{i\phi}\cdot\zeta)|_{\phi=0}\;,\quad \hat S\in
\ccal^\infty(X_P^c)\;.
\end{equation}

\begin{prop} For  $1\le j\le m$,
\begin{itemize}

\item[\rm (i)] $\ \Xi_j:\hcal^2_N(X_P^c) \to \hcal^2_N(X_P^c)$;

\item[\rm (ii)] The lifted monomials $\hat\chi_{\alpha}^P \in \hcal^2_N(X_P^c)$ satisfy
$\Xi_j \hat\chi_{\alpha}^P = \alpha_j \hat\chi_{\alpha}^P$
($\al\in NP$).

\end{itemize}
\end{prop}

\begin{proof}Item (i) follows from the fact that the $\T$ action is holomorphic and
commutes with $\frac{\d}{\d\theta}$. For the case $N=1$, (ii)
follows immediately from (\ref{mhat}) and (\ref{actionX}).  For
$\al\in NP$, $N>1$, we write $\al=\be^1\cdots\be^N$ with $\be^k\in
P$, and the conclusion then follows from the first case and the
product rule.\end{proof}

Furthermore, we recall that
\begin{equation}\label{dtheta} \frac{\d}{\d\theta} :\hcal^2_N(X_P^c) \to
\hcal^2_N(X_P^c)\;,\qquad \frac{1}{i}\frac{\d}{\d\theta}\hat s_N =
N\hat s_N \quad \mbox{for }\ \hat s_N\in
\hcal^2_N(X_P^c)\;.\end{equation}
Thus, we have the joint eigenvalue problem \eqref{jointeigen} with the commuting operators:
\begin{equation}
\hat{I}_{j}=\Xi_{j},\;\;j=1,\ldots,m,\;\;
\hat{I}_{m+1}=\frac{p}{i}\frac{\partial}{\partial \theta}
-\sum_{j=1}^{m}\Xi_{j}.
\label{Qtorus}
\end{equation}
The monomials $\wh{\chi}_{\wh{\alpha}}$ are the joint eigenfunctions
of $(\hat{I}_{1},\ldots,\hat{I}_{m+1})$
for the joint eigenvalues $\wh{\alpha} \in \Lambda_{P}$.

\begin{rem}  The vector fields $\Xi_j$ can be constructed geometrically as follows (see
\cite{Gu}): Let $\xi_j=\frac{\partial}{\partial\phi_j}$ ($1\le
j\le m$) denote the Hamiltonian vector fields generating the $\T$
action on $M_P^c$.  There is a natural contact 1-form $\al$ on
$X_P^c$ determined by the Hermitian connection;  a key property of
$\alpha $ is that $d \alpha = \pi^* \omega$  (see \cite{Z}). We
use $\al$ to define the horizontal lifts of the Hamilton vector
fields $\xi_j$:
$$\pi_* \xi^h_{j} = \xi_j,\;\;\; \alpha(\xi^h_j) = 0.$$
The vector fields $\Xi_j$ are then given by: \begin{equation*}
\Xi_j = \xi^h_j + 2 \pi i \langle \mu_P^c\circ\pi, \xi_j^* \rangle
\frac{\partial}{\partial \theta} =\xi^h_j + 2\pi i
(\mu_P^c\circ\pi)_j\, \frac{\partial}{\partial \theta}.
\end{equation*} \end{rem}
(Here, $\xi_j^* \in \R^m$ is the element of the Lie algebra of
$\T$ which acts as  $\xi_j$ on $M_P$.)

\subsection{\label{FATV} Fourier and Toeplitz analysis   on toric varieties}

Under the lifted torus action of ${\bf T}^m$ on $X_P$
generated by the $\Xi_j$, we can further decompose $\lcal^2(X_P^c) $
and $\hcal^2(X_P^c)$ into representations of ${\bf T}^m$. Since
the action commutes with $S^1$, we combine the two actions into an
action of ${\bf T}^{m + 1}$ on $X_P$.

The ${\bf T}^{m + 1}$ action on $\hcal^2$ is multiplicity free and we
have
$$\hcal^2(X_P^c) = \bigoplus_{\hat{\alpha} \in \Lambda_{P}} \C
\hat{\phi}_{\hat{\alpha}}.$$ We will also need to decompose the
action on $X_P$. To this end, we introduce the anti-Hardy space
$\overline{\hcal}^2(X_P^c)$ of anti-CR functions, i.e. solutions
of $\partial_b f = 0$. Of course, a Hilbert basis is given by the
complex-conjugate monomials $\bar{\hat{\chi}}_{\hat{\alpha}}.$

Because the ${\bf T}^{m + 1}$ action  is generated by vector
fields, products of  eigenfunctions are also eigenfunctions.
Hence, the orthonormal mixed monomials
$$\hat{\chi}_{\hat{\alpha}, \hat{\beta}}(x)
=\hat{\chi}_{\hat{\alpha}}\bar{\hat{\chi}}_{\hat{\beta}}$$ are
eigenfunctions of eigenvalue $\hat{\alpha} - \hat{\beta}$ for
$\{\hat{I}_{1},\ldots,\hat{I}_{m+1}\}$ defined in \eqref{Qtorus}.
We now claim that
these mixed monomials furnish a Hilbert basis of $\lcal^2(X_P^c).$

\begin{prop}  We have: $\lcal^2(X_P^c) = \bigoplus_{\hat{\alpha}, \hat{\beta} \in
\Lambda_{P}} \C \hat{\chi}_{\hat{\alpha}, \hat{\beta}}. $ \end{prop}

\begin{proof} It suffices to show that the closure of the algebra
generated by the $\hat{\chi}_{\alpha, \beta}$ equals the space
$C(X_P^c)$ of continuous functions on $X_P$. Since the span of the
monomials is a *-algebra,  it suffices  by Stone-Weierstrass to
show that such polynomials separate points. But this is clear
since $L_P^c$ is ample:  for any $z \in M_P$ there exists $N, s_N
\in H^0(M_P^c, L_P^{c N})$ such that $|s_N(z)| \not= 0.$ It
follows that $|\hat{s}_N(x)| \not=0$ for any $x$ over $z$.

\end{proof}

It follows that the joint spectrum of $(\hat{I}_{1},\ldots,\hat{I}_{m+1})$ is given by
$$Spec (\hat{I}_{1},\ldots,I_{m+1}) = \Lambda_{P} - \Lambda_{P} = \Z^{m +
1}.$$ However, the multiplicity of a lattice point $\hat\ga \in
\Z^{m+1}$ is infinite since any $(\hat{\alpha}, \hat{\beta})$ with $\hat{\alpha} -
\hat{\beta} = \hat\ga$ corresponds to an eigenvector. This reflects the fact
that $(\Xi_1, \dots, \Xi_m, \frac{\partial}{\partial \theta})$ is
not an elliptic system.

\subsection{Fourier multipliers and Toeplitz Fourier multipliers}

We now define the analogues of convolution operators or Fourier
multipliers:

\begin{defin} An operator $\Pi F \Pi $ on $\hcal^2(X_P)$ will be called a
{\it Fourier multiplier} if it satisfies the following
(equivalent) conditions:

\begin{itemize}

\item $F$ may be expressed as a function $F(D)$ of the commuting
system of operators $D = (\Xi_1, \dots, \Xi_m,
\frac{\partial}{\partial \theta})$.

\item Its eigenfunctions are the monomials
$\hat{\chi}_{\hat{\alpha}, \hat{\beta}}$.

\end{itemize}

\end{defin}

As emphasized above, our main interest is in holomorphic
functions. The relevant definition is:

\begin{defin} An operator $\Pi F \Pi $ on $\hcal^2(X_P)$ will be called a
{\it Toeplitz Fourier multiplier} if its eigenfunctions are the
holomorphic  monomials $\hat{\chi}^{P}_{\hat\alpha }$.

\end{defin}

It is obvious that $\Pi F \Pi$ is a Toeplitz Fourier multiplier if
$F$ is a Fourier multiplier, but the converse is obviously not
necessarily true. In fact, it is sometimes not obvious  how to
extend a Toeplitz Fourier multiplier to all of $\lcal^2(X)$ and to
characterize it there.

\section{Pull back of the projective space \szego kernel}\label{s-szego}

Theorem \ref{M} gives a relation between the \szego kernel of a
toric variety and the very simple kernel obtained by pulling back
the projective space \szego kernel under a holomorphic embedding.
As mentioned in the introduction, the pull back kernels are the
key objects in the proofs of all of our results.  In this section,
we study the pull back kernel on any positive line bundle over any
\kahler manifold, and express it  as a modification of the \szego
kernel. In the next section, we will show that on a toric variety
this relation can be inverted to give a construction of the \szego
kernel.

Consider any polarized algebraic manifold $(M, L)$, and assume $L$
is very ample. Choose a basis $\{S_0, \dots, S_{d}\}$ of $H^0(M, L)$ and let $\Phi$
denote the associated embedding into projective space. That is, we write $S_j = f_j
e_L$ relative to a local frame $e_L$ and put $\Phi(z) = [f_0(z), \dots,
f_{d}(z)]$.  Recalling that $L=\Phi^*\ocal(1)$, we can equip $L$
with the metric $\Phi^*h_\FS$.  We also give $M$ the metric
$\omega = \Phi^* \omega_{\FS}$, and we let $\Pi_N$ denote the
orthogonal projection onto $H^0(M, L^N)$ with respect to these
metrics. Also, let $\Pi= \sum_{N=1}^{\infty} \Pi_N$ denote the
\szego kernel.

 The embedding $\Phi$ determines the lift $\iota:X \to S^{2d +1}$.
Recalling (\ref{iota}), we let
\begin{equation} \Pi_{\iota} (x, y) = \Pi(\iota(x), \iota(y)) = \sum
\Pi_N^{\CP^{d}} (\iota(x),\iota(y))
= \frac{d!}{\left(1- \langle \iota(x), \overline{\iota(y)}
\rangle\right)^{d+1}},\quad x,y\in X\; .
\end{equation}
denote the pullback to $X$ of the \szego kernel of
$\ocal(1)\to \CP^{d}.$ The question we study in this section is the relation
of this simple kernel to the \szego kernel $\Pi=\Pi^M$ with respect to $\iota^*
\omega_{\FS}.$ The following proposition shows that $\Pi_{\iota}$ is a Toeplitz
operator. We assume basic knowledge of such operators and their symbols (see
e.g.  \cite{Gu2}). We denote by  $S^{k}(M, \N)$ the space of
semiclassical symbols $s: M \times \N \to \C$ of the form
$$ s_N(z) \sim N^{ k}
\sum_{j = 0}^{\infty} s_{-j}(z)  N^{-j}.$$

We let \begin{equation}
\wt \Pi_1^{M}(x,y):=\frac{1}{d+1}
\Pi_1^{\CP^{d}}(\iota(x),\iota(y))=
\langle\iota(x),\overline{\iota(y)}\rangle\;.\label{Pi1}\end{equation}
For simplicity of notation, we shall write
\begin{equation}\Pi_1^N(x,y)=\big[\wt
\Pi_1^{M}(x,y)\big]^N.\label{Pi1N}\end{equation}  We note that
$\Pi_1^N(x,x)=1.$

\begin{prop} \label{PAM} Let $(M, L)$ be a polarized algebraic manifold as
above. Then there exists a semi-classical symbol $a_N \in
S^{-m}(M, \N)$ with principal symbol $s_{-m} = 1$ so that
$$   \Pi_1^N = \Pi_N a_N \Pi_N  + R_N\quad \mbox{where }\
\|R_N\|_{\rm HS} = O(N^{-k})\;\;\; \forall k\;. $$ Here,
$\|R_N\|_{\rm HS}^2 = \mbox{\rm Trace}\, R_N^* R_N.$
\end{prop}

This proposition will be used in the proof of Theorem \ref{M}.
The main step in the proof expresses $\Pi$ as the composition of
$\Pi_{\iota}$ and of a Toeplitz operator.

\begin{lem}\label{PIOTA}  Let $\Pi$ denote the \szego projector
associated to the metric $\iota^* \omega_{\FS}$.  Then, there
exist $A \in \Psi^m(X)$ with $[A, D_{\theta}] = 0$ such that $\Pi
\sim  \Pi  A  \Pi \Pi_{\iota} =  \Pi  A \Pi_{\iota}  $ modulo
smoothing operators $\Pi R \Pi.$
\end{lem}

  Here, $\Psi^m(X)$ denotes the space of pseudodifferential operators of
order $m$ on $X$. An operator of the form  $\Pi A \Pi$ where $A
\in \Psi^m(X)$ for some $m$ is a Toeplitz operator in the sense of
\cite{BG}. A smoothing Toeplitz operator is a smoothing operator
of the form $\Pi R \Pi$. It follows then that $\Pi_N = T_N
\Pi_1^N$ where $T_N$ is a Toeplitz operator.

 \subsubsection{The projective \szego kernels}

Since our formula for the \szego kernel of a general toric variety
involves the pullback of the projective space \szego kernel, we
recall how the latter is defined.
  We note that
\begin{equation*}\wh\chi_\al(x) =  x^{\hat\al^p}
\;,\end{equation*} and hence the \szego kernel $\Pi_p^{\CP^m}$ for
the orthogonal projection is given by:
\begin{equation}\label{szego-proj}\Pi_{p}^{\CP^m}(x,y)=\sum_{|\al|\le
p}\frac{1}{\|\chi_\al\|_{\CP^m}}\wh \chi_\al(x) \overline {\wh
\chi_\al(y)} =\frac{(p+m)!}{p!}\sum_{|\al|\le p} {p\choose \al}
x^{\hat \al^p} \bar y^{\hat\be^p} =\frac{(p+m)!}{p!}\langle x,\bar
y\rangle^p \;,\end{equation} for $x,y\in S^{2m+1}$. (The sum
$\sum_{p=0}^\infty \Pi_{p}^{\CP^m}$ is the usual \szego kernel for
the sphere.)

Note  that when $M_\Si=M_{p\Si}=\CP^m$,  we have the circle
bundles $X_{p\Si}\to \CP^m$, for $p\ge 1$. When $p = 1$, $X_\Si =
S^{2m + 1}$ while for $p>1$, it is the lens space $X_{p\Si} =
S^{2m + 1}/ \{e^{2 \pi i/p}\}$. The latter statement follows from
the fact that homogeneous polynomials of degree $p$ are well
defined on (and separate points of) the quotient by the cyclic
group of $p$-th roots of unity.

 \subsubsection{Boutet de Monvel -Sj\"ostrand parametrix} To prove
 Lemma \ref{PIOTA}, we need to recall some background on
 parametrices for $\Pi$.
It was proved by Boutet de Monvel and Sj\"ostrand \cite{BS} (see
also the Appendix to \cite{BG})  that $\Pi$ is a complex Fourier
integral operator of positive type,
\begin{equation}\label{FIO}  \Pi \in I_c^0(X \times X, {\mathcal C}) \end{equation}
associated to a positive canonical relation ${\mathcal C}$. For
definitions and notation concerning complex FIO's we refer to
\cite{MS, BS, BG}. The real points of ${\mathcal C}$ form the
diagonal $\Delta_{\Sigma \times \Sigma}$ in the square  of the
symplectic cone
\begin{equation} \label{SYMCONE} \Sigma = \{r \alpha_x: r > 0, x
\in X\},
\end{equation}
where $\alpha$ is the connection form.
 We refer to \cite{BG} (see Lemma 4.5 of the Appendix).
   Moreover, in \cite{BS}  a parametrix is constructed for $\Pi,$
from which it follows (see \cite{Z}) that
\begin{equation} \label{oscint}\Pi_N(x,y) \sim N \int_0^{\infty}
\int_0^{2\pi} e^{  N ( -i\theta + t  \psi( r_{\theta} x,y))}
s(r_{\theta} x,y, Nt )\, d\theta\, dt \,,\end{equation}  where
$s(x,y,t ) \sim \sum_{k = 0}^{\infty} t ^{m -k} s_k(x,y)\in
S^m(X\times X\times \R_+)$ is  a classical symbol of order $m$.
Here, `$\sim$' means modulo a rapidly decaying term (i.e., a term
whose $\ccal^j$ norms are $O(N^{-k})$ for all $j,k$).

To describe the phase in (\ref{oscint}), we let $e_L$ be a
nonvanishing holomorphic section of $L$ over an open $U\subset M$,
and  consider the  analytic extension' $a(z,w)$ of $a(z,z):=a(z)=
\|e_L(z)\|_h^{-2}$  in $U\times U$ such that
$a(w,z)=\overline{a(z,w)}$ on $U\times U$. Using coordinates
$(z,\theta)$ for the point $x=e^{i\theta} a(z)^\half e_L(z)\in X$,
we have
\begin{equation}\psi(x_1, x_2) = -1 + e^{i (\theta_1 - \theta_2)}
\frac{a(z_1, {z}_2)}{\sqrt{a(z_1)}
\sqrt{a(z_2)}}\;.\label{psi-general}\end{equation}

Again assuming that  the metric  $\om$ on $M$ is the pull-back of
$\om_\FS$, we claim  that  the phase equals:
\begin{equation}  \psi_(x,y)= -1 +
\langle
\iota(x),\overline{\iota(y)}\rangle\;,\label{psiX}\end{equation}
where $\iota:X\to \C^{2d+1}$ is the lift of $\Phi$ given by
$$\iota(z,\theta)=e^{i\theta}
\left(\sum|f_j(z)|^2\right)^{-\half}\big(f_0(z),\dots,f_d(z)\big)\;.$$

To see this, let us recall the \szego kernel of the hyperplane
section bundle $\ocal(1)\to \CP^d$ over projective space with the
Fubini-Study metric. We take $U=\{z_0\ne 0\}\approx \C^d$, and we
consider the local frame $e =z_0$. Using the local coordinates
$[1,z_1,\dots,z_m]\mapsto (z_1,\dots,z_m)\in\C^d$, we then have
$a^{\CP^d}(z)=\|e(z)\|^{-2}=1+\sum_{j=1}^d |z_j|^2$. Hence
$a^{\CP^d}$ has the real-analytic extension
\begin{equation} \label{a-projective} a^{\CP^d}(z,w)=1+\sum_{j=1}^d z_j\bar
w_j\;,\end{equation} and (\ref{psiX})  follows.

\begin{rem} In the case of a
toric variety $M_P^c$, (\ref{psiX}) becomes
$$\psi(x,y) = -1 +\sum_{\al\in
P}\wh m_\al^P(x)\overline{\wh m_\al^P(y)}\;. $$
\end{rem}

\subsection{Proof of Lemma \ref{PIOTA}}\label{pf-TOEP}

We assume as known  the following  facts about complex Fourier
integral operators and Toeplitz operators:
\begin{itemize}

\item $A \in I^m_c(X \times X, {\mathcal C})$ possesses a principal symbol $\sigma_A$ which is a
half-density (times a Maslov factor) along the underlying
canonical relation ${\mathcal C}$. (In the Toeplitz case, it is a
symplectic spinor.)

\item  We can compose operators in $I^m_c(X \times X, {\mathcal C})$ on the
left and right by elements $B \in \Psi^k(X)$ and $\sigma_{AB} =
\sigma_A \sigma_B = \sigma_{BA}$. (The same is true of Toeplitz
operators.)

\item If $A \in I^m_c(X \times X, {\mathcal C})$ and if $\sigma_A = 0$, then
$A \in I^{m-1}_c(X \times X, {\mathcal C})$ (and also for Toeplitz
operators).

\item $\Pi$ and $\Pi_{\iota}$ are elliptic in that their symbols
are nowhere vanishing. See \cite{BS} for the symbol in the complex
FIO sense and \cite[\S 11]{BG} for the symbol in the Toeplitz
sense.

\end{itemize}

We now give the proof:

\begin{proof} As mentioned above (\ref{FIO}),   $\Pi$ is a complex
Fourier integral operator associated to a positive canonical
relation $C$, whose real points form the isotropic relation
$\Delta_{\Sigma \times \Sigma} \subset T^*X \times T^*X$ (the
diagonal).
 We  observe that also $\Pi_{\iota} \in
I^*(X \times X, {\mathcal C})$. This follows immediately from the
fact, show by (\ref{psiX}), that $\Pi$ and $\Pi_{\iota}$ are
complex Fourier integral distributions with precisely the same
phase functions. Since the underlying canonical relation  is
parametrized by the phase, they both belong to the same class of
Fourier integral distributions.

Now, the principal symbol $\sigma_{\Pi}$ of $\Pi$, viewed as a
complex Fourier integral distribution, is a nowhere vanishing
1/2-density on ${\mathcal C}$ which is computed in
\cite[Prop.~4.8]{BS}. Alternatively, viewed as a Toeplitz operator
in the sense of \cite{BG} (see Chapter 11), its symbol is an
idempotent symplectic spinor .
 Similarly, the principal symbol $\sigma_{\Pi_{\iota}}$ of $\Pi_{\iota}$
is the pull back under $\iota$ of the nowhere vanishing  symbol of
$\Pi^{\CP^d}$.

By our normalization, $\Pi_{\iota}$ has order $- m$ (since its
amplitude is a constant independent of $N$). We therefore begin by
seeking $A_0 \in \Psi^m(X)$ such that $[A_0, D_{\theta}] = 0$ and
such that  $\Pi- \Pi A_0 \Pi_{\iota}$ is of order $-1$. We first
find $a_0 \in C^{\infty}(M)$ such that $\sigma_{\Pi} = a_0
\sigma_{\Pi_{\iota}}$ and choose $A_0$ so that $ [A_0, D_{\theta}]
= 0$ and so that  $\sigma_{A_0} = a_0.$ Existence of such an $A_0$
follows by ellipticity of $\Pi_{\iota}$ and by averaging; see also
\cite[Prop.~2.13]{BG}.   Thus,  the principal symbol of order $0$
of $\Pi - \Pi A_0 \Pi_{\iota}$ equals zero, i.e. $\Pi - \Pi A_0
\Pi_{\iota} \in I^{-1}(X \times X, {\mathcal C}).$ We denote its
principal symbol by $\sigma_{-1}.$ We then seek $A_{-1} \in
\Psi^{m}(X)$ so that $ [A_{-1}, D_{\theta}] = 0,$ and so that $\Pi
- \Pi A_0 \Pi_{\iota} - \Pi A_{-1} D_{\theta}^{-1} \Pi_{\iota} \in
I^{-2}(X \times X, {\mathcal C}).$ Here, we note that $\Pi
D_{\theta} \Pi$ is an elliptic Toeplitz operator; since $\Pi
\Pi_{\iota} = \Pi_{\iota}$, the expressions $ D_{\theta}^{-1}
\Pi_{\iota}$ are well-defined.
 It suffices to choose $a_{-1} = \sigma_{A_{-1}} \in
C^{\infty}(M)$ so that $a_{-1} \sigma_{\Pi_{\iota}} =
\sigma_{-1}.$ We continue in this way to obtain $a_{-j} \in
C^{\infty}(M) $  always using that $\sigma_{\Pi_{\iota}}$ is
nowhere vanishing. By a Borel summation argument, we can find $A
\in \Psi^{m}(X)$ with the above commutation properties so that
$\Pi A \Pi - \sum_{j = 0}^M \Pi A_{-j} D_{\theta}^{-j} \Pi$ is a
Toeplitz operator of order $-M - 1.$ Then
$$ \Pi - \Pi A \Pi_{\iota} \in I^{-\infty}(X \times X, {\mathcal C}).
$$
\end{proof}

We remark that an  alternative to the observation that the two
complex FIO's have the same phase function is that, by the  choice
of the \kahler form $\iota^*\omega_{\FS}$,   $\iota$ is a
symplectic (as well as holomorphic)
 embedding $M
\subset \CP^d$. Hence, the  pull back operator $\iota^*$ carries
the class of FIO's in the class of $\Pi^{\CP^d}$ to those in the
class of $\Pi$.

\subsection{Proof of Proposition \ref{PAM}}

To prove the proposition, we need to relate the complex FIO $\Pi A
\Pi_{\iota} $ to semi-classical Toeplitz operators. We will define
the necessary terms as the proof proceeds.

We will need the following asymptotic formula for symbols. In the
language of Berezin-Toeplitz operators, it amounts to computing
the Berezin transform between covariant and contravariant symbols
of a Toeplitz operator.

\begin{lem} \label{BEREZIN} Let $\sigma_N \sim N^{ k} \sum_{j = 0}^{\infty} s_{-j} N^{-j}$ be  a   semiclassical  symbol of order $k$.  Then
there exists a complete asymptotic expansion
$$ N^{-m}(\Pi_N\sigma_N\Pi_N)(z,z)=  \sum_{j = 0}^{l-1} b_{-j}(z) N^{k-j} +
r_N^l(z) \qquad (l\ge 1)\;,$$ where $b_0 = s_0,b_{-1} = \Delta s_0
+ s_{-1}, \dots$, and in general where $ b_{-j}$ is a sum of
differential operators applied to $s_0, s_{-1}, \dots, s_{-j}.$
Also,
 $ \|D_z^n r_N^l (z)\| \leq
C_{nl} N^{k - l}$. \end{lem}

\begin{proof} Apply the method of stationary phase using
(\ref{oscint})--(\ref{psi-general}) exactly as in the proof of
\cite[Theorem~1]{Z} (which is the case $\sigma=1$).\end{proof}

\subsubsection{Proof of Proposition \ref{PAM}}

 By Lemma \ref{PIOTA},  there exists $A \in \Psi^{-m}(X)$
such that  $\Pi_{\iota} = \Pi A \Pi.$ Since $[A, D_{\theta}] = 0$,
it follows  that there exists a symbol $a_N$ with $\|\Pi_1^N-
\Pi_N a_N \Pi_N\|_{\rm HS}=O(N^{-\infty})$.    We may determine
$a_N$ by using Lemma \ref{BEREZIN}. Indeed, we have
\begin{equation}\label{diag} \Pi_1^N(z, z) = 1  \sim \Pi_N a_N \Pi_N
(z,z)\end{equation} modulo functions $r_N$ which tend to zero
rapidly in $\ccal^k(M)$. It follows that $a_0 = 1$ (and then the
rest of the coefficients may be determined recursively, e.g. $
a_{-1} = 0, a_{-2} = - D_2 a_0$, and so on).
\qed

\subsubsection{Alternate proof of Proposition \ref{PAM}}
One could avoid using the calculus of complex Fourier integral
operators or Toeplitz operators in the proof of Proposition
\ref{TOEP} by further developing the calculus of semi-classical
Toeplitz operators directly from the Boutet de Monvel- Sj\"ostrand
parametrix as follows:

We first simplify the expression (\ref{oscint}) by using the
complex method of stationary phase  to eliminate the integrals in
the parametrix. The critical point set of  the phase $\Phi
(\theta, \lambda; x, y):= -i \theta + \lambda \psi(r_{\theta} x,
y)$ is given by
\begin{equation} \left\{ \begin{array}{l}d_{\lambda} \Phi (\theta, \lambda; x, y) =
\psi(r_{\theta} x, y) = 0 \iff  e^{i \theta} \langle \iota(x), \iota( y) \rangle= 1,\\
\\  d_{\theta} \Phi (\theta, \lambda; x, y) = -i + d_{\theta} \lambda \psi(r_{\theta}
x, y) = 0 \iff  \lambda e^{i \theta} \langle \iota(x), \iota( y)
\rangle = 1.
\end{array} \right. \end{equation} It is easy to see by the Schwartz inequality that a
real critical point exists if and only if $x = y,$ in which case
$\theta = 0, \lambda = 1$ and we obtain the familiar expansion
along the diagonal. When $x \not= y$ we deform the  contour to
$|\zeta| = e^{\tau}$ ($\zeta = e^{i \theta + \tau}$)  so that  $
e^{i \theta + \tau} \langle \iota(x), \iota( y) \rangle = 1$. This
is possible as long as $\langle \iota(x), \iota( y) \rangle \not=
0$, as happens near the diagonal, where the parametrix is valid.
Because the phase is linear in $\lambda$ it is clear that the
critical point is  non-degenerate if and only if $\langle
\iota(x), \iota( y) \rangle  \not= 0$ and that the Hessian
determinant equals $|\langle \iota(x), \iota( y) \rangle |^2.$ On
the critical set the phase equals $- i \theta + \tau = \log
\langle \iota(x), \iota( y) \rangle,$ hence  we have
\begin{equation} \label{oscintN} \Pi_N(x, y) =
e^{N \log \langle \iota(x), \iota( y) \rangle} S_N(x, y) +
W_N(x,y), \end{equation} where $S_N(x,y) \sim \sum_{k =
0}^{\infty} N^{m - k} S_k(x,y),$ and where $W_N(x,y)$ is a smooth
uniformly rapidly decaying function. Here, we have absorbed the
remainder in the parametrix construction as well as the remainder
in the stationary phase expansion of the parametrix in $W_N$. Note
that the first term may be smaller than the second outside a
tubular  neighborhood of radius $N^{-1/2}$ of the diagonal.

As above, we use Lemma \ref{BEREZIN} to find a symbol $a_N$ so
that (\ref{diag}) holds and hence  the kernels  $\Pi_1^N$ and
$\Pi_N a_N \Pi_N$  agree on the diagonal modulo smoothing symbols.
Next, we note that both $\Pi_1^N(x,y)$ and $\Pi_N a_N \Pi_N(x,y)$
are complex oscillatory functions with  common phase
\begin{equation} \label{PHASE} \Psi(z, w) =  \log \langle \iota(x), \iota(y) \rangle.
\end{equation} In the case of $\Pi_1^N$, this  follows from (\ref{Pi1}). Indeed, we
simply have:
\begin{equation} \label{PiN} \Pi_1^N(x,y) = e^{N  \log \langle \iota(x), \iota(y)
\rangle }. \end{equation} In the case of $\Pi_N a_N \Pi_N$, we
apply the method of complex  stationary phase to the integral
formula
\begin{equation} \Pi_N a_N \Pi_N(x,y) \sim \int_{M_P} e^{N \Psi(u; x,y)} S_N(x,u) a_N(u)
S_N(u, y) dV(u), \end{equation}
 coming from (\ref{oscintN}),  with $$ \Psi(u; x,y) = \langle \iota(x), \iota( u)
\rangle + \langle \iota(u), \iota( y) \rangle.$$ It follows that
 there exists an amplitude
$A_N$ defined near the diagonal such that
\begin{equation}\label{oscintS} \Pi_N a_N \Pi_N(x, y) = A_N(x,y) e^{N\log \langle
\iota(x), \iota( u) \rangle} + V_N(x,y),\end{equation} where $V_N$
is a new smoothing operator.

Recalling that $X\subset L^*$, we extend $A_N$ to $L^*\times L^*$
so that it is of the form
$$A_N(z,\la; w, \la') = (\la\bar\la')^N \tilde A_N(z,w)\;,$$
where we use a local holomorphic frame to write $x = (z, \la),\ y
= (w, \la')\in L^*$.  We note that $\tilde A_N(z,w)$ is
holomorphic in $z$ and anti-holomorphic in $w$ near the diagonal.
To see this, we first conclude from the construction in \cite{BS}
using the $*$-product that the symbol $s(x,y,t) \sim \sum_{k =
0}^{\infty} t ^{m -k} s_k(x,y)$ in (\ref{oscint}) extends to a
symbol on $L^*\times L^*$ that  is holomorphic in $x$ and
anti-holomorphic in $y$.  It follows by the stationary phase
method described above that the same is true for the symbol
$S_N(x,y)$ in (\ref{oscintN}), and then by (\ref{oscintS}) that
the same is also true for $A_N$ as claimed.  (Note that in terms
of coordinates $(z,\theta)$ on $X$, the function $z\mapsto
A_N(z,\theta;w,\theta')
=e^{iN(\theta-\theta')}\sqrt{a(z)a(w)}\,\tilde A_N(z,w)$ is not
holomorphic in $z$.)

The amplitude of $\Pi_1^N - \Pi_N a_N \Pi_N$ becomes
$(\la\bar\la')^N (1 - \tilde{A}_N(z,w))$, where we can write
$1-\tilde A_N \sim B_0 +B_1 N^{-1}+B_2 N^{-2}+\cdots$.  Since
$\Pi_1^N - \Pi_N a_N \Pi_N(z,z)=O(N^{-\infty})$ by our choice of
$a_N$ above, it follows that $B_j(z,z)=0$ (for all $j$).  Since
$B_j(z,w)$ vanishes on the diagonal and is holomorphic in $z$ and
anti-holomorphic in $w$, it must be identically 0.  Hence,
$\Pi_1^N - \Pi_N a_N \Pi_N(x, y)$ is a smoothing operator.\qed

\section{\label{SKTV} \szego kernels on toric varieties}

The purpose of this section is to give a special construction of
the \szego kernels
\begin{equation} \label{Sztor} \Pi_N^{M_P}(x,y)  = \sum_{\alpha \in NP}
\frac 1{\|\chi_\al^P\|_{M_P}^2}\; \chi_\al^P(x) \overline{
\chi_{\alpha}^P(y)} \;, \end{equation}
 of a toric variety. The  first step, Lemma
\ref{PQ}, is to  give an exact formula for $\Pi_N^{M_P}$ as the
composition of a certain {\it Toeplitz-Fourier multiplier} denoted
$({\mathcal P}{\mathcal Q})^{-1}$ and the pull back of the
Fubini-Study kernel under a monomial embedding. The pulled-back
Fubini-Study kernel is very simple to analyze and is the raison
d'etre of our method. The Toeplitz-Fourier multiplier requires
more work. In the second step, Lemma \ref{TOEP}, we prove that
this operator is a ${\bf T}^{m+1}$-invariant  {\it Toeplitz
operator} of order $m $ on $X_P$ (modulo a smoothing operator).

\subsection{Proof of Theorem \ref{M}: The exact formula}\label{s-exact}

Our exact formula for the \szego kernel of a toric variety
involves two ingredients: The first is the kernel $\Pi_1^N$ given by
(\ref{Pi1})--(\ref{Pi1N}).
The second ingredient is the Fourier multiplier ${\mathcal M}$
defined by the eigenvalues (\ref{FM}). We lift the multiplier to
$X$ as follows:  We observe that each sequence of  functions can
be re-defined as a single function on the homogenized lattice cone
$\Lambda_P=\bigcup_{N = 1}^{\infty} \wh{NP}$ defined in Section \ref{FAnalysis}.
We define $\pcal,\qcal:\La_P\to \R_+$ by
\begin{eqnarray*}{\mathcal P}(\wh\alpha^{N}):= {\mathcal
P}_N(\alpha) &=& \# \{(\beta_1,
\dots, \beta_N): \beta_j \in P, \beta_1 + \cdots + \beta_N =
\alpha \},\\ {\mathcal Q}(\wh\alpha^{N}) :=  {\mathcal
Q}_N(\alpha) &=& \int_{X_P} |\wh\chi^P_\al (x)|^2 d\vol_{X_P}(x)\;,\end{eqnarray*}
and  define Toeplitz Fourier multipliers:

\begin{equation}\begin{array}{l} {\mathcal P}(D) \wh \chi_\al^P = {\mathcal P}(\wh\alpha^N)
\wh \chi_\al^P,\;\;\; {\mathcal Q}(D) \wh \chi_\al^P = {\mathcal
Q} (\wh\alpha^N) \wh \chi_\al^P,\;\;\; \alpha \in N P \\ \\
{\mathcal M}:= ({\mathcal P} {\mathcal Q})^{-1}. \end{array}
\end{equation}
(Recall that $\pcal_N(\al)\ge 1$ by (\ref{solid}), and hence $\pcal\qcal(\wh
\al_N)>0$.) We write ${\mathcal M}_N$ for the restriction of ${\mathcal M}$ to
functions in the $N$-th subspace.

We now specialize the kernels (\ref{Pi1})--(\ref{Pi1N}) to the case $M=M_P$:
\begin{equation}
\wt \Pi_1^{M_P}(x,y)=
\langle\iota_P(x),\overline{\iota_P(y)}\rangle = \sum_{\alpha \in
P} \wh m_\al^P(x) \overline{\wh
m_\al^P(y)}\;,\qquad \Pi_1^N(x,y)=\big[\wt
\Pi_1^{M_P}(x,y)\big]^N\label{Pi1MP}\;,\end{equation}
where
$\iota_P :X_P \to S^{2d+1}$ is the lift of the monomial
embedding  $M_P \hookrightarrow \CP^d$, as described in \S
\ref{s-lifting}. We shall choose the constants $c_\al = 1$ for all
$\al\in P$, so that $\wh m_\al^P = \wh\chi_\al^P$.  (Except for
Lemma~\ref{PQ} below, all our results hold without change for
arbitrary $\{c_\al\}$.)

We now obtain the exact formula for the \szego kernel  in Theorem
\ref{M}. The proof is very simple.

\begin{lem}\label{PQ} The \szego kernels of a toric variety $M_P$
factor as follows:
$$\Pi_N^{M_P} = {\mathcal M}_N \circ \Pi_1^N. $$ \end{lem}

\begin{proof}

First, we have by definition, $$\Pi_N^{M_P}(x, y ) = \sum_{\alpha
\in N P} \frac{1}{{\mathcal Q}_N(\alpha)}  \wh\chi_\al^P(x)
\overline{\wh\chi_\al^P(y)}\;.$$

On the other hand, by definition of the partition function, we
also have
$$ \Pi_1^N (x,y) = \sum_{\alpha \in N P} {\mathcal
P}_N(\alpha)  \wh\chi_\al^P(x) \overline{\wh\chi_\al^P(y)}\;.
$$
We note that the $N$-th power of $\Pi_1^{M_P}$ gives the sum
over the correct set of exponents but does not have the correct
normalizing coefficients. We need to divide each term by
${\mathcal P}_N(\alpha) {\mathcal Q}_N(\alpha)$ to adjust the
coefficients. That is just what the lemma claims.
\end{proof}

Of course,
$$({\mathcal P} {\mathcal Q})^{-1} \Pi_1^N = \Pi_N^{M_P} ({\mathcal P} {\mathcal Q})^{-1} \Pi_N^{M_P} \Pi_1^N,$$
so ${\mathcal M}$ need only be defined on ${\mathcal H}^2$.

\subsection{Proof of Theorem \ref{M}: ${\mathcal M}$ is a Toeplitz operator }

The explicit formula of Lemma \ref{PQ} is difficult to use
without a detailed analysis of the multiplier ${\mathcal M}$. In
this section, we prove that ${\mathcal M}$ is a Toeplitz operator,
i.e. possesses a symbol. This has numerous implication for
${\mathcal M}$, in particular, that it has polyhomogeneous
expansions along rays.

We begin with some   background on symbols. We consider symbols
$\sigma$  of the form  $\sigma(z, D_{\theta})$, where $\sigma(z,
N)$  is a semiclassical symbol.
 Here,  we say  $\sigma \in S^k(M
\times \N)$ is a semiclassical symbol of order $k$ if
$$\sigma_N(z)=\sigma(z,N) = N^k \sum_{j = 0}^{l-1} a_j(z) N^{-j} +
r_N^l(z),\;\; \mbox{with}\;\; \|D_z^n r_N^l (z)\| \leq C_{nl} N^{k
- l}\qquad (l\ge 1), $$ where $a_j\in\ccal^\infty(M)$. It is a
{\it smoothing symbol\/} if $\|D_z^n \sigma_N (z)\| \leq
C_{n}N^{-l} $ for all $n,l$. The {\it Toeplitz operator\/}
associated to a symbol $\sigma$ is the operator $\Pi {\sigma}(z,
D_\theta)\Pi $, where  $D_{\theta}$ denotes the symbol of
$\frac{\partial}{\partial \theta}$ (the $S^1$ generator).  Its
symbol is the polyhomogeneous function on the symplectic  cone
$\Sigma = \{(x, r \alpha_x) : r > 0\} \subset T^*X$ given by
$$\sigma(z, p_{\theta})  \sim \sum_{j = 0}^{\infty} a_j(z)
p_{\theta}^{k-j} \;. $$ Since $\sigma$ commutes with the $S^1$
action, we have $\Pi \sigma \Pi = \sum_N \Pi_N \sigma_N \Pi_N$.

The
goal of this section is to prove that $\Pi_N (\pcal\qcal )^{-1}
\Pi_N $ is a semi-classical Toeplitz operator in the following
sense:

\begin{lem}\label{TOEP} There exists a ${\bf T}^m$-invariant symbol $\sigma(z, D_{\theta})$
so that ${\mathcal M} = \Pi \sigma(z, D_{\theta}) \Pi$.
Equivalently,  there exists a symbol $\sigma_N$ of order $m$ with
principal symbol equal to $1$ and a smoothing operator $R_N$ so
that
$$\Pi_N (\pcal\qcal )^{-1} \Pi_N = \Pi_N \sigma_N \Pi_N + R_N\;.$$ \end{lem}

\begin{proof}

 In the toric case, $\Pi_N (\pcal\qcal )^{-1} \Pi_N$ is the inverse
of $\Pi_N (\pcal\qcal ) \Pi_N = \Pi_1^N$ on $\hcal^2(X_P)$. From
Lemma \ref{PAM} we may write $\Pi_N (\pcal\qcal ) \Pi_N  \sim
\Pi_N a_N \Pi_N$ modulo smoothing operators.

We note that we may invert $\Pi a \Pi$ in the class of Toeplitz
operators; i.e.,  there exists a symbol $\sigma_N$ such that
\begin{equation}\label{INVERSE}  \Pi_N \sigma_N \Pi_N \circ \Pi_N a_N \Pi_N \sim
\Pi_N \end{equation} modulo smoothing operators. Such an inverse
symbol exists since $a_0 = 1$. The algebraic formalism in which
the inverse is calculated is that of $*$-products of semiclassical
symbols. We recall that composition of Toeplitz operators defines
a $*$-product on semiclassical symbols by the formula
\begin{equation} \Pi_N a_N \Pi_N \circ \Pi_N b_N \Pi_N \sim \Pi_N
a_N * b_N \Pi_N. \end{equation} The formula for $*$ may be worked
out directly from the parametrix (\ref{oscint}) and the inverse
can be computed from this formula (see \cite{Gu2, BS}).

By (\ref{INVERSE}) we obtain a symbol $\sigma_N$ (with principal
symbol equal to $1$) such that
$$\Pi_N \sigma_N \Pi_N\Pi_N (\pcal\qcal ) \Pi_N \sim \Pi_N\;.$$
Multiplying both sides by $\Pi_N (\pcal\qcal )^{-1} \Pi_N$, we
conclude the proof. \end{proof}

\begin{rem}
We emphasize that the distinguishing features of  the toric case
in Proposition  \ref{PAM} and Lemma  \ref{TOEP} are the exact
factorization and the fact that the operator $\Pi ({\mathcal
P}{\mathcal Q})\Pi$ mediating between $\Pi$ and $\Pi_{\iota}$ is
invertible. This is due to the fact that sections of $L$ generate
the ring $\oplus_{N=1}^{\infty} H^0(M, L^N)$ in the toric case. In
general they do not and the exact representation $\Pi = B
\Pi_{\iota}$ in the toric case  is only valid modulo smoothing
operators. \end{rem}

\subsection{Polytope characters}

As a corollary of  Lemma
\ref{TOEP}, we obtain the following formula for the
character of the torus action:

\begin{prop}\label{aha}  There exists a symbol $\sigma_N \in S^{m}(M, \N)$
with principal symbol equal to $1$ and a smoothing operator $R_N$
such that for $t\in \T $, we have
\begin{eqnarray*}\chi_{NP}(t) &=& \int_{X_P} \sigma_N(x) \Pi_1^{N}( t
\cdot x,x) \,d\vol_{X_P} (x)\\&&+ \int_{X_P} \int_{X_P} R_N(y,x)
\Pi_1^{N}(t \cdot x,  y) \,d\vol_{X_P} (x)\,d\vol_{X_P} (y).
\end{eqnarray*}
\end{prop}

\begin{proof} We shall use
the elementary formula for the polytope character:
\begin{equation}\chi_{NP} (e^{i\phi}) = \int_{M_P} \Pi_N(e^{i\phi}\cdot  x, x )
dV(x)\;. \label{elem}\end{equation}
Applying
$\Pi_1^N=\Pi_N (\pcal\qcal ) \Pi_N$ to the identity of Lemma
\ref{TOEP}, we obtain
$$\Pi_N  = (\Pi_N \sigma_N \Pi_N + R_N) \Pi_1^N. $$
Now let $T_{t}$ denote the translation operator $f(x)\mapsto
f(t \cdot x)$ on $\lcal^2(X_P)$. Since $[T_t, \Pi_N] = 0$, we
then have
$$T_t\Pi_N= T_t\Pi_N(\sigma_N+R_N)\Pi_1^N =
\Pi_NT_t(\sigma_N+R_N)\Pi_1^N\;.$$ Therefore,
\begin{eqnarray*} &&\int_{X_P} \Pi_N( t\cdot x,
x) \,d\vol_{X_P} (x)\ =\ \mbox{Trace\,}   T_t \Pi_N \ =\
\mbox{Trace\,}\Pi_NT_t(\sigma_N+R_N)\Pi_1^N\\ &&\qquad\qquad=\
\mbox{Trace\,}T_t\sigma_N\Pi_1^N +
\mbox{Trace\,}T_tR_N\Pi_1^N\\&&\qquad\qquad=\ \int_{X_P}
\sigma_N(t\cdot x) \Pi_1^{N}(t\cdot x,x) \,d\vol_{X_P} (x)
\\&&\qquad\qquad\qquad+ \int_{X_P}\int_{X_P} R_N(t\cdot x,y)
\Pi_1^{N}(y,x) \,d\vol_{X_P} (x) \,d\vol_{X_P} (y)\;.
\end{eqnarray*}
Making the change of variables $x\mapsto t\inv\cdot
x$, we then obtain the desired formula from  (\ref{elem}).\end{proof}

\end{document}